\documentclass[11pt]{amsart}
\usepackage{latexsym}
\usepackage{amsfonts}

\usepackage{amsmath,amsthm,amssymb}

\usepackage{epsf,psfrag}

\newtheorem{theor}{Theorem}

\newtheorem{theore}{Theorem}

\newtheorem{thm}{Theorem}[section] \newtheorem{lem}[thm]{Lemma}
\newtheorem{cor}[thm]{Corollary} 
\newtheorem{prop}[thm]{Proposition} \theoremstyle{definition}
\newtheorem{defn}[thm]{Definition}
\newtheorem{notation}[thm]{Notation}
 \newtheorem{rem}[thm]{Remark}
\newtheorem{exmp}[thm]{Example}

 \newtheorem{prob}[thm]{Problem}

\begin{document}

\title[Translation equivalence in free groups]{Translation equivalence
  in free groups}

\author[I.~Kapovich]{Ilya Kapovich}

\address{\tt Department of Mathematics, University of Illinois at
  Urbana-Champaign, 1409 West Green Street, Urbana, IL 61801, USA
  \newline http://www.math.uiuc.edu/\~{}kapovich/} \email{\tt
  kapovich@math.uiuc.edu}

\thanks{The first author acknowledges the support of the Max Planck
  Institute of Mathematics in Bonn. The first and the third author
  were supported by the NSF grant DMS\#0404991 and the NSA grant
  DMA\#H98230-04-1-0115. The fourth author was supported by the NSF
  grant DMS\#0405105}

\author[G.~Levitt]{Gilbert Levitt} \address{\tt Laboratoire de
  Mathematiques Nicolas Oresme, CNRS UMR 6139, Universite de Caen, BP
  5186, 14032 Caen Cedex, France} \email{\tt
  Gilbert.Levitt@math.unicaen.fr}

\author[P.~Schupp]{Paul Schupp}

\address{\tt Department of Mathematics, University of Illinois at
  Urbana-Champaign, 1409 West Green Street, Urbana, IL 61801, USA
  \newline http://www.math.uiuc.edu/People/schupp.html}
\email{schupp@math.uiuc.edu}

\author[V.~Shpilrain]{Vladimir Shpilrain} \address{\tt Department of
  Mathematics, The City College of New York, New York, NY 10031, USA
  \newline http://www.sci.ccny.cuny.edu/\~{}shpil} \email{\tt
  shpil@groups.sci.ccny.cuny.edu}

\begin{abstract}
  Motivated by the work of Leininger on hyperbolic equivalence of
  homotopy classes of closed curves on surfaces, we investigate a
  similar phenomenon for free groups. Namely, we study the situation
  when two elements $g,h$ in a free group $F$ have the property that
  for every free isometric action of $F$ on an $\mathbb{R}$-tree $X$
  the translation lengths of $g$ and $h$ on $X$ are equal.
\end{abstract}


\maketitle

\section{Introduction}\label{intro}

Let $S$ be a closed oriented surface of negative Euler characteristic
and let $\gamma$ be a free homotopy class of essential closed curves
on $S$. If $\rho$ is a hyperbolic metric on $S$ then $\gamma$ contains
a unique curve $c$ of minimal $\rho$-length. We denote this length by
$\ell_{\rho}(\gamma)$. The curve $c$ is a closed geodesic on $(S,
\rho)$ and $\ell_{\rho}(\gamma)$ is the translation length of any
representative of $\gamma$ in the action corresponding to $\rho$ of
$G=\pi_1(S)$ on $\mathbb{H}^2 = \tilde{S}$. There is an obvious
identification between the set of nontrivial conjugacy classes $C$ of
$G$ and the set of free homotopy classes of essential closed curves on
$S$ and we shall not distinguish between the two.  Thus each marked
hyperbolic structure on $S$ defines a so-called \textit{marked length
  spectrum} $l:C \to \mathbb{R}$. It is well-known and easy to see
that a marked hyperbolic structure on $S$, considered as a point of
the Teichm\"uller space of $S$, is uniquely determined by its marked
length spectrum.

The dual situation, however, is different.  For $\gamma_1, \gamma_2
\in C$ we say that $\gamma_1$ is \emph{hyperbolically equivalent to}
$\gamma_2$, denoted $\gamma_1 \equiv_h \gamma_2$, if for every
hyperbolic structure $\rho$ on $S$ we have $\ell_{\rho}(\gamma_1)=
\ell_{\rho}(\gamma_2)$.  In more algebraic terms, for two conjugacy
classes $\gamma_1,\gamma_2 \in C$ we have $\gamma_1 \equiv_h \gamma_2$
if for every discrete and co-compact isometric action of $G$ on
$\mathbb{H}^2$ the translation lengths of $\gamma_1$ and $\gamma_2$
are equal. It can happen that $\gamma_1 \ne \gamma_2^{\pm 1}$ and yet
$\gamma_1 \equiv_h \gamma_2$. The main source of hyperbolic
equivalence comes from ``trace identities" in $SL(2, \mathbb{C})$.  A
number of interesting new results about hyperbolic equivalence were
recently obtained by Chris Leininger~\cite{Le}.

We are interested in investigating a similar phenomenon for free
groups. In this context the Teichm\"uller space is replaced by the
Culler-Vogtmann outer space~\cite{CV}, so that instead of actions on
$\mathbb{H}^2$ we consider free and discrete actions on
$\mathbb{R}$-trees. Recall that if $G$ is a group acting by isometries
on an $\mathbb{R}$-tree $X$ and $g\in G$ then the \emph{translation
  length} $\ell_X(g)$ is defined as:
\[
\ell_X(g) = \inf_{x \in X} d (x,gx) .
\]
It is easy to see that $\ell_X(g)$ only depends on the conjugacy class
of $g$ and that in the definition above the infimum can be replaced by
a minimum. Thus $\ell_X(g)=0$ if and only if $g$ fixes a point of $X$.
This discussion naturally leads us to the following definition:

\begin{defn}
  Let $F$ be a finitely generated free group and let $g, h \in F$ be
  elements of $F$. We say that $g$ and $h$ are \emph{translation
    equivalent in $F$}, denoted $g \equiv_{t} h$, if for every free
  and discrete isometric action of $F$ on an $\mathbb{R}$-tree $X$ we
  have $\ell_X(g) = \ell_X(h)$.
\end{defn}

It is obvious that $\equiv_t$ is an equivalence relation on $F$.
Applying the above definition to the case where $X$ is the Cayley
graph of $F$ with respect to some free basis of $F$ implies that every
$\equiv_t$-equivalence class is the union of finitely many conjugacy
classes in $F$.  Clearly, if $g \equiv_{t} h$ in $F$ and $\phi:F \to
F_1$ is an injective homomorphism to a free group $F_1$ then $\phi(g)
\equiv_{t} \phi(h)$ in $F_1$.  Indeed, suppose $F_1$ acts freely and
discretely by isometries on an $\mathbb R$-tree $X$. Then, by
restriction, we get a free and discrete action of $\phi(F)$ on $X$,
and via a twist by $\phi$, we also get a free and discrete action of
$F$ itself on $X$. Namely, $f\cdot p:=\phi(f) p$, where $p\in X, f\in
F$.  Since $u$ and $v$ are translation equivalent in $F$, it follows
that $\ell_X(u)=\ell_X(v)$, that is,
$\ell_X(\phi(u))=\ell_X(\phi(v))$.

A phenomenon related to but different from translation equivalence was
studied by Smillie and Vogtmann~\cite{SV} who, given an arbitrary
finite set of conjugacy classes in a free group, constructed
multi-parametric families of free discrete actions where the
translation length of each conjugacy class from this set remains
constant through the family. Results similar in spirit to those of
\cite{SV} were also obtained by Cohen, Lustig and Steiner~\cite{CLS}.

The notion of translation equivalence is also related to the space of
geodesic currents on a free group. In~\cite{Ka04} the first author
studies the properties of an \emph{intersection form} 
\[
I: FLen(F)\times Curr(F)\to \mathbb R.
\] 
Here $FLen(F)$ is the
(non-projectivized) space of hyperbolic length functions corresponding
to free and discrete isometric actions of a $F$ on $\mathbb R$-trees
and $Curr(F)$ is the space of \emph{geodesic currents} on $F$, that is
the space of $F$-invariant positive Borel measures on the set of all
pairs $(x,y)$, where $x,y\in \partial F$ and $x\ne y$. Similarly to
Bonahon's notion~\cite{Bo} of the intersection number between geodesic
currents on hyperbolic surfaces, it turns out that if $\eta_g$ is the
``counting'' current corresponding to a nontrivial $g\in F$ and if
$\ell\in FLen(F)$ is a length function then $I(\ell,\eta_g)=\ell(g)$.
Thus $g\equiv_t h$ in $F$ if and only if for every $\ell\in FLen(F)$
we have $I(\ell,\eta_g)=I(\ell,\eta_h)$. Therefore the notion of
translation equivalence, in a sense, measures the degeneracy of the
intersection form $I$ with respect to its second argument. We refer
the reader to~\cite{Ka04} for a detailed discussion on this topic.

The first natural problem is to demonstrate that there are nontrivial
instances of translation equivalence in free groups. We provide two
different sources of translation equivalence: one based on trace
identities in $SL(2,\mathbb C)$ and another based on power
redistribution for certain products of translation equivalent
elements. Both methods can be iterated and used to produce arbitrarily
large finite collections of distinct conjugacy classes in $F$ that are
pairwise translation equivalent. Both of these sources can also be
used to produce hyperbolic equivalence in the context of
$\mathbb{H}^2$-actions, although some distinctions do arise as will be
pointed out later. (See Example~\ref{distinct}.)

Another natural question is to give a more algebraic and combinatorial
characterization of translation equivalence.

Recall that an isometric action of a group $G$ on an $\mathbb{R}$-tree
$X$ is \emph{very small}~\cite{CL} if the following conditions hold:
\begin{enumerate}
\item The action is \emph{small}, that is, arc stabilizers do not
  contain free subgroups of rank two.
  
\item Stabilizers of tripods are trivial.
  
\item For any $g \in G$ and for each $n\ne 0$ the fixed sets of $g$
  and $g^n$ are equal.
\end{enumerate}
In particular, every free action and, more generally, an action with
trivial arc stabilizers, is very small. Results of
Cohen-Lustig~\cite{CL} and Bestvina-Feighn~\cite{BF93} imply that an
action of $F_n$ on an $\mathbb{R}$-tree is very small if and only if
this action represents a point in the standard length functions
compactification of the outer space. Thus a very small action can
always be approximated in the sense of length functions by a sequence
of free simplicial actions.

\begin{notation} If $A$ is a basis of a free group $F$ any element  $g\in F$
  is represented by a unique reduced word $w_g$ over the alphabet
  $A^{\pm 1}$.  The \emph{length of $g$ with respect to the basis
    $A$}, denoted $|g|_A$, is the number of letters in $w_g$.  A word
  is \emph{cyclically reduced} if all its cyclic permutations are
  reduced.  Any reduced word $w$ can be uniquely factored as $w = c u
  c^{-1} $ where $u$ is cyclically reduced. If $w_g = c u c^{-1}$ is
  such a factorization, then $|u|_A$ is called \emph{the cyclically
    reduced length of $g$ with respect to $A$} is denoted by
  $||g||_A$. Note that $||g||_A=\ell_X(g)$ where $X$ is the Cayley
  graph of $F$ with respect to $A$.  A \emph{cyclic word} in $A^{\pm
    1}$ is the set of all cyclic permutations of a cyclically reduced
  word.  There is a canonical identification between the set of cyclic
  words in $A^{\pm 1}$ and the set of conjugacy classes in $F$.
  
  If $w$ is a cyclic word consisting of all cyclic permutations of a
  nontrivial word $u$, and if $v$ is a word in $A^{\pm 1}$, we define
  the \emph{number of occurrences of $v$ in $w$} as the number of
  those $i$, $0\le i<|u|$, such that the infinite word $uuu\dots$
  begins with $u_iv$, where $u_i$ is the initial segment of $u$ of
  length $i$.  If $w$ is a cyclic word and $x,y \in A^{\pm 1}$ we use
  $n_A(w;x,y)$ (or just $n(w;x,y)$ if $A$ is fixed) to denote the
  total number of occurrences of the subwords $x y$ and $y^{-1}
  x^{-1}$ in $w$. Thus $n(w;x,y) = n(w;y^{- 1},x^{- 1}) =n(w^{-
    1};x,y)$.  Similarly, in this case if $w$ is a cyclic word and
  $x\in A^{\pm 1}$, we denote by $n(w;x)$ the total number of
  occurrences of $x$ and $x^{-1}$ in $w$. Thus again
  $n(w;x)=n(w;x^{-1}) = n(w^{-1};x)$. If $[g]$ is a nontrivial
  conjugacy class in $F$ and $w$ is the unique cyclic word over $A$
  representing $[g]$, we denote $n_A([g];x):=n(w;x)$, where $x\in
  A^{\pm 1}$.
\end{notation}

In studying automorphisms of free groups, Whitehead automorphisms and
the Whitehead graph of a cyclically reduced word play a major role.
See Lyndon-Schupp \cite{LS} for a detailed discussion. Note that in
\cite{LS} Whitehead graphs are called \emph{star graphs}.

\begin{defn}[Whitehead graph]
  Let $w$ be a nontrivial cyclic word in $F(A)$. The \emph{Whitehead
    graph $\mathcal{W}_A(w)$ of $w$ with respect to $A$} is the
  labelled undirected graph defined as follows. The vertex set of
  $\mathcal{W}_A(w)$ is $A^{\pm 1}$. If $x, y \in A^{\pm 1}$ are such
  that $x \ne y$, there is an edge in $\mathcal{W}_A(w)$ between $x$
  to $y$ with label $n(w;x,y^{-1})$.

  If $[g]$ is a nontrivial conjugacy class in $F$, then $[g]$ is
  represented by a unique cyclic word $w$ in $F(A)$. Then the
  \textit{Whitehead graph $\mathcal{W}_A ([g])$ of $[g]$ with respect
    to $A$} is defined as $\mathcal{W}_A(w)$.  Note that
  $\mathcal{W}_A (w) = \mathcal{W}_A (w^{- 1})$ for any nontrivial
  cyclic word $w$.

\end{defn}

We obtain the following result:

\begin{theor}\label{A}
  Let $F$ be a finitely generated free group and let $g, h \in F$ be
  nontrivial elements.
  
  Then the following statements are equivalent:
\begin{enumerate}
\item $\ell_X(g) = \ell_X(h)$ for every very small action of $F$ on an
  $\mathbb{R}$-tree $X$.

\item $\ell_X(g) = \ell_X(h)$ for every free action of $F$ on an
  $\mathbb{R}$-tree $X$.

\item $g \equiv_{t} h$ in $F$.
  
\item $||g||_A = ||h||_A$ for every free basis $A$ of $F$.
  
\item $\mathcal{W}_A ([g]) = \mathcal{W}_A ([h])$ for every free basis
  $A$ of $F$.  That is, the conjugacy classes $[g]$ and $[h]$ have the
  same Whitehead graphs with respect to $A$.
  
\item $n_A([g];x)=n_A([h];x)$ for every free basis $A$ of $F$ and for
  every $x\in A$.
\end{enumerate}
\end{theor}

Theorem~\ref{A} immediately yields a more combinatorial version of
translation equivalence:

\begin{cor}
  Let $F$ be a finitely generated free group with a free basis $A$ and
  let $g, h \in F$. Then the following conditions are equivalent:
\begin{enumerate}
\item $g \equiv_{t} h$ in $F$.
  
\item $||\phi(g)||_A =||\phi(h)||_A$ for every automorphism $\phi$ of
  $F$.
  
\item $||\phi(g)||_A = ||\phi(h)||_A$ for every injective endomorphism
  $\phi$ of $F$.
  
\item $||\phi(g)||_B =||\phi(h)||_B$ for every free group $F_1$ with a
  free basis $B$ and for every injective homomorphism $\phi : F \to
  F_1$.

\end{enumerate}
\end{cor}

The a priori weakest condition in Theorem~\ref{A} is condition (4)
which deserves further comment. Let $S$ be a closed surface as before.
If $\gamma$ and $\delta$ are free homotopy classes of closed curves on
$S$, we use $i(\gamma, \delta)$ to denote the \emph{geometric
  intersection number} of $\gamma$ and $\delta$. We say that free
homotopy classes $\gamma_1, \gamma_2$ of closed curves on $S$ are
\emph{simple intersection equivalent} if
$i(\gamma_1,[c])=i(\gamma_2,[c])$ for every essential \textit{simple
  closed curve} $c$ on $S$.  It is easy to see that hyperbolic
equivalence implies simple intersection equivalence.  Surprisingly
however, the converse is not true as was recently proved by Chris
Leininger~\cite{Le}.

It is therefore natural to consider the analogue of simple
intersection equivalence for free groups. If $G = \pi_1(S)$ and $c$ is
an essential simple closed curve on $S$, then $c$ defines a splitting
of $G$ as either an amalgamated product or as an HNN-extension over a
cyclic subgroup. Let $X$ be the Bass-Serre tree corresponding to that
splitting. It is easy to see that for any $g \in G$ we have $i([g],
[c])= \ell_X(g)$. The difficulty is that for a free group $F$, a
single element of $F$, even if it is a primitive one, does not define
a splitting of $F$.  A free basis $A$ of $F$ does, however, define
such a splitting. Namely, the splitting of $F$ as the multiple
HNN-extension of the trivial group with stable letters corresponding
to the elements of $A$ and the Bass-Serre tree $X$ of this splitting
is precisely the Cayley graph of $F$ with respect to $A$. Then for any
$g \in F$ we have $\ell_X(g) =||g||_A$, the cyclically reduced length
of $g$ with respect to $A$.  For a free basis $A$ of $F$ and an
element $g \in F$ we can therefore define the \textit{intersection
  number} $i(g,A):=||g||_A$.  By analogy with the surface group case
we say that $g_1,g_2 \in F$ are \emph{simple intersection equivalent}
in $F$ if for every free basis $A$ of $F$ we have $i(g_1,A) =
i(g_2,A)$. Unlike in the surface group case, Theorem~\ref{A} says that
in free groups simple intersection equivalence is the same as
translation equivalence.

We can now state the main sources of translation equivalence in free
groups that we have discovered so far.

If $F$ is a finitely generated free group and $u,v\in F$ are
nontrivial elements, we say that $u$ and $v$ are \emph{trace
  equivalent} or \emph{character equivalent} in $F$, denoted $u
\equiv_{c} v$, if for every representation $\alpha: F\to SL(2,\mathbb
C)$ we have $tr(\alpha(u))=tr(\alpha(v))$. Character equivalent words
come from the so-called ``trace identities'' in $SL(2,\mathbb C)$ and
are quite plentiful (see, for example, \cite{Ho}). A corollary of
Theorem~\ref{A} together with a result of Horowitz~\cite{Ho} is:

\begin{theor}\label{B}
  Let $F$ be a finitely generated free group and suppose that
  $u\equiv_c v$ in $F$. Then $u\equiv_t v$ in $F$.
\end{theor}

A particularly interesting example of character equivalence comes from
two-variable ``palindromic reversing''.  Let $w(x,y) \in F(x,y)$ be a
freely reduced word. We denote by $w^R(x,y)$ the word $w(x,y)$ read
backwards, but without inverting the letters. Thus $w^R(x,y) = (w(x^{-
  1}, y^{- 1}))^{- 1}$. We prove:

\begin{theor}\label{C}
  Let $F$ be a free group of rank $k \ge 2$ and let $w(x,y) \in
  F(x,y)$ be a freely reduced word. Then for any $g, h \in F$ we have
\[
w(g,h) \equiv_{t} w^R(g,h) \quad \text{ in } F. 
\]
\end{theor}

We shall give a direct proof of the above statement as well as a proof
via character equivalence.

A very different source of translation equivalence is given by:

\begin{theor}\label{D}
  Let $F$ be a free group of rank $k \ge 2$ and let $g,h \in F$ be
  such that $g \equiv_{t} h$ but $g \ne h^{-1}$.  Then for any
  positive integers $p, q, i, j$ such that $p+q = i+j$ we have
\[
g^p h^q \equiv_{t} g^i h^j \quad \text{ in } F.
\]
\end{theor}
Theorem~\ref{B} states that $\equiv_c$ implies $\equiv_t$. However, it
turns out that Theorem~\ref{D} does not hold for character
equivalence, as demonstrated by Example~\ref{distinct} below. This
example shows that $\equiv_t$ does not imply $ \equiv_c$ and that,
although character equivalence and translation equivalence are closely
related phenomena, they are not the same and translation equivalence
is more general.

The following is an analogue of a result of Randol~\cite{Ra} about
hyperbolic surfaces.

\begin{cor}\label{cor:C}
  Let $F$ be a free group of rank $k \ge 2$. Then for any integer $M
  \ge 1$ there exist elements $g_1, \dots, g_M \in F$ such that $g_i
  \equiv_{t} g_j$ in $F$ and such that for $i \ne j$ $g_i$ is not
  conjugate to $g_j^{\pm 1}$.
\end{cor}

\begin{proof}
  Let $M \ge 1$ be an integer and let $(a, b, \dots )$ be a free basis
  of $F$. Consider $g=a$ and $h = b a^{-1} b^{-1}$. Put $g_i = g^i
  h^{2M-i} = a^i b a^{i-2M} b^{-1}$ for $i = 1, \dots, M$.  Then by
  Theorem~\ref{D} $g_i \equiv_{t} g_j$ in $F$. On the other hand $g_i$
  is not conjugate to $g_j^{\pm 1}$ for $i \ne j$. This can be seen,
  for example, by observing that $g_i$ and $g_j^{\pm 1}$ have distinct
  images in the abelianization of $F$.
\end{proof}

\begin{rem}
  Theorem~\ref{C} can be iterated to produce other examples with the
  same properties as in Corollary~\ref{cor:C}. Namely, let
  $\phi_i:F(x,y) \to F(x,y)$ be injective endomorphisms of $F(x,y)$
  for $i = 1, \dots, N$.  Let $\phi := \phi_N \circ \dots \circ
  \phi_1$. For $i = 1, \dots, N-1$ put $\psi_i : = \phi_N \circ \dots
  \circ \phi_{i+1}$ and $\theta_i = \phi_i \circ \dots \circ \phi_1$.
  Then $\phi = \psi_i \circ \theta_i$.  Let $\psi_i (x) = u_i(x,y)$,
  $\psi_i(y)=v_i(x,y)$, $\theta_i(x)=r_i(x,y)$ and $\theta_i(y)
  =s_i(x,y)$.  Let $w(x,y) = \phi(x)$. Then $w = \psi_i(\theta_i(x)) =
  \psi_i(r_i(x,y)) = r_i(u_i(x,y),v_i(x,y))$.  Let $w_i =
  r_i^R(u_i(x,y), v_i(x,y))$. Theorem~\ref{C} implies that, $w(g,h)
  \equiv_{t} w_i(g,h)$ for any $g, h \in F$ for $i =1, \dots, N-1$. It
  is possible, for each $N \ge 1$, to choose the endomorphisms
  $\phi_i$ and then elements $g, h \in F$ so that $w_1, \dots,
  w_{N-1}$ are pairwise non-conjugate in $F$.
  
  Corollary~\ref{cor:C} also follows from Theorem~\ref{B} and the
  result of Horowitz~\cite{Ho} establishing (via a more complicated
  family of words) a similar result for character equivalence.
  Moreover, as we observe later in Corollary~\ref{cor:C'}, it is
  possible to generalize the proof of Corollary~\ref{cor:C} to many
  non-free groups.
\end{rem}

Using $SL_2$ trace identities, we show that Theorem~\ref{C} and a
version of Theorem~\ref{D} also hold for standard hyperbolic
equivalence. Via a limiting argument we conclude that these statements
also apply to the tree actions that occur in the Thurston boundary of
the Teichm\"uller space of a closed hyperbolic surface, either
orientable or non-orientable. Recall that each point $\mu$ in the
Thurston boundary of the Teichm\"uller space is a measured lamination.
There is an $\mathbb{R}$-tree $X_{\mu}$, dual to the lift of this
lamination, that comes equipped with a small isometric action of the
fundamental group of the surface. In the case of a non-orientable
surface not all such actions are very small.

\begin{theor}\label{E} Let $S$ be a possibly non-orientable closed surface  of negative Euler
  characteristic and let $G =\pi_1(S)$. Then the following hold:
\begin{enumerate}
\item For any $g, h \in G$ and for any $w(x,y) \in F(x,y)$ and for any
  tree action $\mu$ of $G$ in the Thurston boundary of the
  Teichm\"uller space of $S$ we have
  \[
  \ell_{X_{\mu}}(w(g,h)) = \ell_{X_{\mu}}(w^R(g,h)).
\]
\item For any conjugate elements $g,h \in G$, for any $p, q > 0$ and
  for any tree action $\mu$ of $G$ in the Thurston boundary of the
  Teichm\"uller space of $S$ we have
  \[
  \ell_{X_{\mu}}( g^p h^q) = \ell_{X_{\mu}}(g^q h^p).
\]
\item If $S$ is orientable then for any conjugate elements $g,h \in
  G$, for each point $\mu\in \overline{\mathcal T(S)}-\mathcal T(S)$
  and for any positive integers $p,q,i,j$ such that $p+q=i+j$ we have
\[
\ell_{X_{\mu}}(g^ph^q)=\ell_{X_{\mu}}(g^ih^j).
\]
\end{enumerate}
\end{theor}

The paper is organized as follows. In Section~2 we discuss Whitehead
graphs and prove Theorem~A. In Section~2 we obtain a direct geometric
proof of Theorem~C about the palindromic sources of translation
equivalence. In Section~4 we use the analysis of possible axis
configurations for compositions of isometries of $\mathbb R$-trees to
provide a geometric proof of Theorem~D. Section~5 contains a
discussion of $SL_2$ trace identities and a proof of Theorem~B. In
Section~6 we analyze how our results apply to tree actions occurring
in the boundary of the Teichm\"uller space and prove Theorem~E. In
Section~7 we discuss various examples and counter-examples, and in
Section~8 we list a number of interesting open problems.

The authors are grateful to the organizers of the ``Geometric Groups
on the Gulf'' Conference in Mobile in February 2004, where the
conversations and discussions eventually leading to the writing of
this paper took place. The authors also thank David Berg, Brian
Bowditch, Enric Ventura and Victor Pan for useful conversations.  The
authors also thank the referee for his careful reading of the paper
and for many helpful comments and suggestions.

\section{Whitehead graphs and a characterization of translation equivalence}

Let $F = F(A) = F(a_1, \dots, a_k)$ be a free group with basis $A = \{
a_1, \dots, a_k \}$ where $k \ge 2$.

Let $x, y \in \{ a_1, \dots, a_k \}^{\pm 1}$ be such that $x \ne
y^{\pm 1}$.  Denote by $\phi_{x, y}$ the Nielsen automorphism of $F$
that sends $x$ to $x y$ and fixes each generator $a_i \ne x^{\pm 1}$.
For simplicity, if $i \ne j, 1 \le i, j \le k$ we use $\phi_{i, j}$ to
denote $ \phi_{a_i, a_j}$.

The following lemma is obvious:

\begin{lem}\label{l1}
  For any cyclic word $w$ and any $x, y \in \{a_1, \dots, a_k \}^{\pm
    1}$ such that $x \ne y^{\pm 1}$ we have
\[
||\phi_{x, y}(w)|| - ||w|| = n(w;x) - 2 n(w;x, y^{- 1}).
\]
\end{lem}

\begin{lem}\label{l2}
  Let $w$ be a cyclic word. Then:
\[
n(w;a_i) = \frac{1}{k} \big( ||w|| + \sum_{j \ne i} (||\phi_{i,
  j}(w)|| - ||\phi_{j, i}(w)||) \big).
\]
\end{lem}

\begin{proof}
  Note that $n(w; x, y^{- 1}) = n(w; y, x^{- 1})$. Therefore

\begin{gather*}
  ||{\phi}_{x,y}(w)||-||w||=n(w;x)-2n(w;x,y^{-1})\\  
  ||{\phi}_{y,x}(w)||-||w||=n(w;y)-2n(w;y,x^{-1})=n(w;y)-2n(w;x,y^{-1}),
\end{gather*}

and so
\[
||\phi_{x, y}(w)|| - ||\phi_{y, x}(w)|| = n(w;x) - n(w; y).
\]

Let $x = a_i$ and $y$ vary over $a_j, j \ne i$. Then summing up the
instances of the above equality for $x = a_i, y = a_j$ we get
\[
(k-1) n(w;a_i) - \sum_{j \ne i} n(w,a_j) = \sum_{j \ne i}
(||\phi_{i, j}(w)|| - ||\phi_{j, i}(w)||)
\]
On the other hand
\[
n(w;a_i) + \sum_{j \ne i} n(w,a_j) = ||w||.
\]
Adding the above formulas we get
\[
k \, n(w;a_i) = ||w|| + 
\sum_{j \ne i}^n (||\phi_{i,j}(w)|| - ||\phi_{j, i}(w)||),
\]
which yields the statement of the lemma.
\end{proof}

\begin{prop}\label{prop:impl}
  Let $F$ be a finitely generated free group and let $u,v \in F$ be
  nontrivial elements such that for every free basis $A$ of $F$ we
  have $||u||_A = ||v||_A$. Then for every free basis $A$ of $F$ the
  Whitehead graphs of $[u]$ and $[v]$ with respect to $A$ are equal.
\end{prop}
\begin{proof}
  Let $A$ be a free basis of $F$. We consider the conjugacy classes
  $[u]$, $[v]$ as cyclic words over $A$. The assumptions of the
  proposition imply that for every automorphism $\phi$ of $F$
\[
|| \phi(u)||_A = ||\phi(v)||_A.
\]

By Lemma~\ref{l2} it follows that for each $x \in A$ we have $n([u];x)
=n([v],x)$. Therefore by Lemma~\ref{l1} for every $x, y \in A^{\pm 1}$
such that $x \ne y^{\pm 1}$ we have
\[
n([u];x, y) = n([v]; x, y ).
\]
For a fixed $x \in A$ and an arbitrary cyclic word $w$
\[
n(w;x) =n(w; x, x) + \sum_{y \ne x^{\pm 1}} n(w; x, y ).
\]
Since $n([ u ]; x) = n([v];x)$ and $n([u]; x, y) = n([v]; x, y)$ for
any $y \ne x^{\pm 1}$, $y \in A^{\pm 1}$, it follows that $n([u]; x,
x) = n([v]; x, x)$.

Thus we have shown that for any $x, y \in A^{\pm 1}$ such that $x \ne
y^{-1}$ we have $n([u] ; x, y) = n ([v]; x, y)$. This means that $[u]$
and $[v]$ have equal Whitehead graphs with respect to $A$, as claimed.
\end{proof}

We can now establish Theorem~\ref{A} from the Introduction:

\begin{theore}
  Let $F$ be a finitely generated free group and let $g, h \in F$ be
  nontrivial elements.
  
  Then the following statements are equivalent:
\begin{enumerate}
\item $\ell_X(g) = \ell_X(h)$ for every very small action of $F$ on an
  $\mathbb{R}$-tree $X$.

\item $\ell_X(g) = \ell_X(h)$ for every free action of $F$ on an
  $\mathbb{R}$-tree $X$.

\item $g \equiv_{t} h$ in $F$.
  
\item $||g||_A = ||h||_A$ for every free basis $A$ of $F$.
  
\item $\mathcal{W}_A ([g]) = \mathcal{W}_A ([h])$ for every free basis
  $A$ of $F$.  That is, the conjugacy classes $[g]$ and $[h]$ have the
  same Whitehead graphs with respect to $A$.
  
\item $n_A([g];x)=n_A([h];x)$ for every free basis $A$ of $F$ and for
  every $x\in A$.
\end{enumerate}
\end{theore}

\begin{proof}
  Let $k$ be the rank of $F$. We may assume that $k\ge 2$ since for
  $k=1$ the statement of the theorem is obvious.
  
  The implications $(1)\Rightarrow (2)$, $(2)\Rightarrow (3)$,
  $(3)\Rightarrow (4)$, $(5) \Rightarrow (6)$ and $(6) \Rightarrow
  (4)$ are obvious. Moreover, $(3)$ implies $(1)$ since by the results
  of \cite{CL} every very small action is the limit (in the sense of
  length functions) of free discrete actions. Thus $(1), (2), (3)$ are
  equivalent.  The implication $(4)\Rightarrow (6)$ follows from
  Lemma~\ref{l2} and the implication $(4)\Rightarrow (5)$ is
  Proposition~\ref{prop:impl}.

  We now show that $(5)$ implies (3). Indeed, suppose that $F$ is
  acting freely, discretely and isometrically on an $\mathbb R$-tree
  $X$. Let $Y=X/F$ be the quotient graph of $X$. Since $X$ is a metric
  tree, the edges $e$ of $Y$ come equipped with the lengths $l(e)$.
  Thus every edge-path in $Y$ has a \emph{length} which is the sum of
  the lengths of the edges of this path. There is an obvious canonical
  identification between $F$ and $\pi_1(Y,y)$ where $y$ is a vertex of
  $Y$.  Choose an orientation $EY=E^+Y\sqcup E^-Y$ on $Y$. Then for
  every maximal tree in $Y$ there is a canonically associated free
  basis of $\pi_1(Y,y)$.

  Given a conjugacy class $[g]$ of $g\in F$, represent it by an
  immersed loop $\gamma$ in $Y$. Then $\ell_X(g)=\sum_{e\in E^+Y}
  n(\gamma; e )l(e)$, where $n(\gamma; e )$ is the number of times
  that $\gamma$ traverses $e$ (in either direction). We prove that (5)
  implies (3) by showing that, given $e\in E^+Y$, the number
  $n(\gamma; e )$ is completely determined by the Whitehead graphs of
  $[g]$ with respect to free bases of $F$.
  
  There are two cases.
  
  First, suppose that $e$ does not separate $Y$. Choose a maximal tree
  not containing $e$, and consider the associated free basis $A$ of
  $F=\pi _1(Y,y)$. Then $n_A(\gamma;e )=n_A([g];a_e)$, where $a_e\in
  A$ is the generator corresponding to $e$.
  
  Suppose next that $e$ separates $Y$, so that $Y=Y_1\cup e \cup Y_2$,
  say with $y \in Y_1$.  Choose a basis $B$ of $F$ associated to any
  maximal tree in $Y$. This basis is partitioned as $B=B_1\sqcup B_2$,
  with $b\in B_i$ if and only if it corresponds to an edge in $Y_i$.
  Then $\displaystyle n(\gamma; e )= \sum_{x\in B_1, y\in
    B_2}n_B([g];x,y)$.

\end{proof}

\section{Palindromic sources of translation equivalence}

Let $F$ be a free group of rank $k \ge 2$ and let $A = \{a_1, \dots,
a_k\}$ be a free basis of $F$. Let $u = u(a_1, \dots, a_k) \in F$ be a
freely reduced word over $A$. We define the \textit{palindromic
  reverse of $u$ with respect to $A$}, denoted $u^R$, as:
\[ u^R:=u(a_1^{- 1}, \dots, a_k^{- 1})^{- 1}. \]
Thus $u^R$ is the word $u$ read backwards without inverting the
letters.

Similarly, we define the palindromic reverse $w^R$ of a cyclic word
$w$ over $A$. Thus $w^R$ is again a cyclic word. Namely, if $w$ is
represented by a cyclically reduced word $u$ then $w^R$ is represented
by the cyclically reduced word $u^R$.

\begin{prop}\label{prop:use}
  Let $F = F (a,b)$ be free of rank two. Then for any cyclic word $w$
  in $F$ over $\{a,b\}^{\pm 1}$ and for any $\phi \in Aut(F)$ we have:
\[
(\phi(w))^R = \phi(w^R).
\]
\end{prop}
\begin{proof}
  For any $\psi \in Aut(F)$ denote by $\overline{\psi}$ the image of
  $\psi$ in $Out(F)$. For a free basis $A =(a,b)$ of $F$ denote by
  $\tau_A$ the automorphism of $F$ defined as $\tau_A(a) = a^{- 1},
  \tau_A(b) = b^{-1}$. It is easy to see that in $Out(F)$ the element
  $\overline{\tau_A}$ commutes with all the elementary Nielsen
  automorphisms with respect to $A$ and hence $\overline{\phi_A}$ is
  central in $Out(F)$. Therefore for any other free basis $B$ of $F$
  we have $\overline{\tau_A} = \overline{\tau_B}$.
  
  Suppose now that $w(a,b)$ is a cyclic word in $F(a,b)$ and $\phi \in
  Aut(F)$. Let $u(a, b) = \phi(w)$. Then by the above observation
  $\phi(w(a^{-1}, b^{-1})) = u(a^{- 1}, b^{- 1})$. Since $w^R =(w(a^{-
    1}, b^{-1}))^{- 1}$ and $u^R =(u(a^{- 1}, b^{- 1}))^{- 1}$, this
  implies that $(\phi(w))^R = \phi(w^R)$.
\end{proof}

\begin{rem}
  It is well-known that $Out(F_2) \cong GL(2, \mathbb{Z})$ and that
  the center of $GL(2, \mathbb{Z})$ is cyclic of order two. Thus in
  fact the outer automorphism $\overline{\tau_A}$ is the only
  nontrivial element of the center of $Out(F_2)$, although we did not
  need this fact in the above proof.
\end{rem}

We can now prove Theorem~\ref{C} from the Introduction:
\setcounter{theore}{2}
\begin{theore}\label{palindr}
  Let $F$ be a free group of rank $k \ge 2$ and let $w(x,y) \in
  F(x,y)$ be a freely reduced word. Then for any $g, h \in F$ we have
\[
w(g,h) \equiv_{t} w^R(g,h) \quad \text{ in } F. \]
\end{theore}

\begin{proof}
  If $g$ and $h$ commute in $F$ then the statement is obvious. Suppose
  now that $g$ and $h$ do not commute and hence $F_1 = \langle g, h
  \rangle$ is a free group of rank two. Let $F$ act freely and
  discretely by isometries on an $\mathbb{R}$-tree $X$. Let $X_1$ be
  the minimal $F_1$-invariant subtree and let $Y = X_1/F_1$ be the
  quotient graph. Then topologically $Y$ is either a wedge of two
  circles or is a $\theta$-graph or $Y$ consists of two disjoint
  circles joined by an edge. In each case $Y$ possesses an involution
  isometry $\sigma$ that leaves a maximal subtree $T$ of $Y$
  invariant, fixes some point $y\in T$ and takes every edge $e$
  outside of a maximal tree to $e^{-1}$.  This is shown in
  Figure~\ref{Fi:inv}.

\begin{figure}[here]
  \epsfxsize=4in \epsfbox{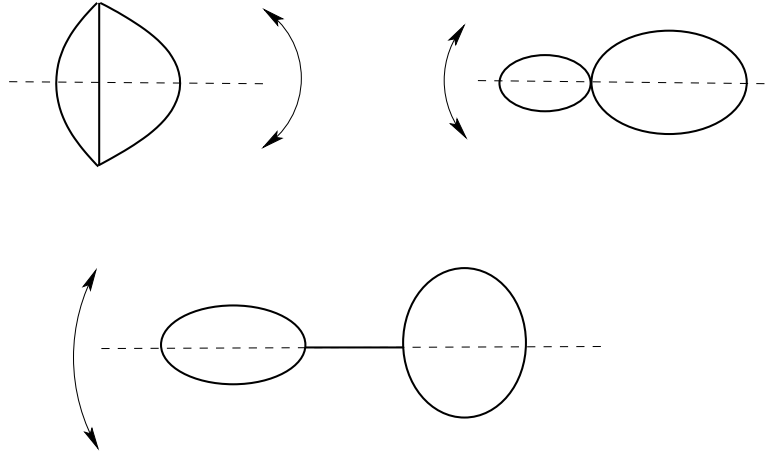}
\caption{``Hyper-elliptic involution"}\label{Fi:inv}
\end{figure}

Thus if $c_1, c_2$ is the basis of $\pi_1(Y,y)$ corresponding to $T$
and $\sigma_{\#}$ denotes the isomorphism of $\pi_1(Y,y)$ induced by
$\sigma$, then for any cyclic word $u(c_1,c_2)$ we have $\sigma_{\#}
u(c_1,c_2) =u(c_1^{-1},c_2^{-1})$. Since $\sigma$ is an isometry of
$Y$, this means that $\ell_X(u(c_1,c_2)) = \ell_X(u(c_1^{- 1}, c_2^{-
  1}))$. On the other hand $u(c_1^{- 1}, c_2^{- 1}) = (u^R(c_1,c_2
))^{-1}$ and so
\[
\ell_X(u( c_1, c_2 )) = \ell_X(u^R(c_1,c_2 )).
\]

The pair $(c_1, c_2)$ is a free basis of $F_1 = F(g,h)$. Write $w(g,h)
=u(c_1, c_2 )$ and $w^R(g,h) = u'(c_1,c_2)$. By
Proposition~\ref{prop:use} $u'= u^R$. Hence
\[
\ell_X(w(g, h)) = \ell_X(u(c_1, c_2)) = \ell_X(u^R(c_1, c_2)) =
\ell_X(w^R(g, h)),
\]
as required.
\end{proof}

\section{Axis diagrams for free actions}

Let $G$ be a group acting by isometries on an $\mathbb{R}$-tree $X$.
Recall that $g\in G$ is called \emph{elliptic} if $\ell_X(g)=0$ and
$g$ is called \emph{hyperbolic} if $\ell_X(g) > 0$. Thus $g$ is
elliptic if and only if it fixes a point of $X$.

For a hyperbolic $g \in G$ put
\[
L_g = \{ x \in X : d(x, g x) = \ell_X(g)\} .
\]
Then $L_g$ is the smallest $g$-invariant subtree of $X$ which is
isometric to a line and on which $g$ acts by a translation of
magnitude $\ell_X( g )$. The set $L_g$ is called the \emph{axis} of
$g$. In this section if an $\mathbb R$-tree $X$ is fixed, we will omit
the subscript and denote the translation length of an isometry $g$ of
$X$ by $\ell(g)$.

The following simple proposition enumerating all the possibilities for
the configuration of $L_{g h}$ with respect to $L_g, L_h$ for two
hyperbolic isometries $g$ and $h$ is essentially a restatement of
Proposition~1.6 of Paulin~\cite{Pau}.

\begin{prop}\label{pau}
  Let $X$ be an $\mathbb R$-tree and let $g, h \in Isom(X)$ be two
  hyperbolic isometries of $X$. Then the following hold:
\begin{enumerate}

\item Suppose that $|L_g \cap L_h| \le 1$ and let $D = d(L_g, L_h)$.
  Then
\[
\ell(g h) = \ell(g) + \ell(h) + 2D.
\]
\item Suppose that $L_g \cap L_h$ is a nondegenerate segment $[x,y]$.

\noindent (a) If the translation directions of $g$ and $h$ on $[x,y]$ coincide then
\[
\ell(g h) = \ell(g) + \ell(h).
\]

\noindent  (b) If the translation directions of $g$ and $h$ on $[x,y]$ are opposite then
\[
\ell(g h) =
\begin{cases}
  \ell(g) + \ell(h) - 2 d(x, y) \quad \text{ if } \ell(g)\ge d(x,y), \ell (h)\ge d(x,y) \\
  |\ell(g) - \ell(h)| \quad \text{otherwise.}
\end{cases}
\]
\item If $L_g$ and $L_h$ are equal or intersect in a ray, then
\[
\ell(gh)=
\begin{cases}
  \ell(g) +\ell(h) \quad \text{ if $g$ and $h$ translate
    in the same direction } \\
  |\ell(g) - \ell(h)| \quad \text{ otherwise.}
\end{cases}
\]
\end{enumerate}
\end{prop}

\begin{rem}
  The case where $L_g\cap L_h$ consists of a single point is omitted
  in Proposition~1.6 of \cite{Pau}. Only the cases where $L_g\cap L_h$
  is empty or contains a nondegenerate segment are explicitly covered
  there. However the proofs for the cases where $L_g\cap L_h$ is empty
  and where $L_g\cap L_h$ is a single point are completely analogous.
\end{rem}

We need another simple fact (see Proposition~1.8 of Paulin~\cite{Pau})

\begin{prop}\label{pau1}
  Let $g,h$ be two elliptic isometries of an $\mathbb R$-tree $X$.
  
  Then
\[
\ell_X(gh)=2 \min \{ d(x,y)| gx=x, hy=y\}.
\]
\end{prop}

The following statement, together with the definition of translation
equivalence, immediately implies Theorem~\ref{D} from the
introduction.

\begin{thm}\label{thm-pq}
  Let $G$ be a group acting isometrically with trivial arc stabilizers
  on an $\mathbb{R}$-tree $X$.  Let $g, h \in G$ be nontrivial
  elements of $G$ such that $g\ne h^{-1}$.

\begin{enumerate}
\item Suppose that $\ell_X(g) =\ell_X(h)>0$.  Then for any positive
  integers $p,q,i,j$ such that $p + q = i + j$ we have
\[
\ell_X(g^p h^q) = \ell_X(g^i h^j).
\]

\item Suppose that $\ell_X(g)=\ell_X(h)=0$.
  
  Then for any integers $p,q, i,j$ such that $g^p,h^q, g^i, h^j$ are
  nontrivial we have
\[
\ell_X(g^ph^q)=\ell_X(g^ih^j)
\]
\end{enumerate}
\end{thm}

\begin{proof}
  
  Part (2) follows directly from Proposition~\ref{pau1}. Indeed, since
  the action of $G$ has trivial arc stabilizers, there are some
  $x,y\in X$ such that $Fix(g)=\{x\}$ and $Fix(h)=\{y\}$. Then for any
  $p,q,i,j$ such that $g^p,g^i,h^q,h^j$ are nontrivial, we have
  $Fix(g^p)=Fix(g^i)=\{x\}$ and $Fix(h^q)=Fix(h^j)=\{y\}$ and hence by
  Proposition~\ref{pau1} $\ell(g^ph^q)=d(x,y)=\ell(g^ih^j)$.
  
  Suppose now that the assumptions of part (1) of Theorem~\ref{thm-pq}
  are satisfied.
  
  Denote $a =\ell(g) = \ell(h)$. If $g = h$ then the statement is
  obvious.  Suppose now that $g \ne h$, so that $g \ne h^{\pm 1}$. Let
  $p,q,i,j \ge 1$ be integers such that $p+q =i+j$.

  Observe that $L_g = L_{g^n}$ and $L_h = L_{h^n}$ for any $n \ge 1$.
  Moreover, in this case $\ell(g^n) = \ell(h^n) = n \ell(h) = n
  \ell(g) = na$.
  
  Suppose first that $L_g \cap L_h$ consists of at most one point. Put
  $D = d(L_g, L_h)$. Then by part (1) of Proposition~\ref{pau}
\begin{gather*}
  \ell(g^p h^q) = \ell(g^p) + \ell(h^q) + 2D = p \ell(g) + q \ell(h)
  + 2 D =\\
  p a + q a + 2 D = (p+q) a + 2 D .
\end{gather*}
Thus we see that $\ell(g^p h^q)$ depends only on $p+q$ and hence
$\ell(g^p h^q) = \ell(g^i h^j)$, as required.

If the intersection of $L_g$ and $L_h$ contains a ray then either $g
h$ or $g h^{- 1}$ fixes a segment of that ray. This is impossible
since the action of $G$ on $X$ has trivial arc stabilizers.

Suppose now that $L_g \cap L_h =[x,y]$ and that $d(x,y) > 0$. If the
translation directions of $g$ and $h$ on $[x,y]$ coincide then by part
(2a) of Proposition~\ref{pau} we have:
\[
\ell(g^p h^q) = \ell(g^p) + \ell(h^q) = (p+q) a = (i+j)a = \ell(g^i
h^j).
\]
Assume now that $g$ and $h$ translate on $[x,y]$ in the opposite
directions.

If $a < d (x,y)$ then $g h$ fixes an arc contained in $[x,y]$,
yielding a contradiction. Hence $\ell(g) = \ell(h)= a \ge d(x,y)$.
Then by part (2b) of Proposition~\ref{pau} we have
\begin{gather*}
  \ell(g^ph^q)=\ell(g^p)+\ell(h^q)-2d(x,y)=\\
  (p+q)a-2d(x,y)=(i+j)a-2d(x,y)=\ell(g^ih^j).
\end{gather*}
\end{proof}

The following lemma is an elementary exercise, but we provide a proof
for completeness.

\begin{lem}\label{ab}
  Let $A$ be a finite abelian group and let $g\in A$ be an element of
  order bigger than four. Then for any $u\in A$ there is an integer
  $i$ such that $g^iu$ has order bigger than four.
\end{lem}
\begin{proof}
  Note that if $a=bc$ in $A$ then the order of any of these three
  elements divides the least common multiple of the orders of the
  other two.
  
  If $u\in \langle g\rangle$, then $u=g^n$ for some $n$ and the
  conclusion of the lemma holds with $i=-n-1$. Suppose now that $u$ is
  not a power of $g$ and that for every $i$ the order of $g^iu$ is at
  most four. Hence for each integer $i$ the order of $g^iu$ is two,
  three or four. If $i$ is an integer then by the observation above
  the orders of $g^iu$ and $g^{i+1}u$ cannot be both even or both
  equal to three. Indeed, in that case the order of $g$ would divide
  either three or four, contrary to the assumption that the order of
  $g$ is bigger than four.  Choose $i$ such that the order of $g^iu$
  is three. Then the order of $g^{i+2}u$ is also three and the orders
  of $g^{i+1}u, g^{i+3}u$ are both even. Hence the order of
  $g^2=g^{i+2}u(g^iu)^{-1}=g^{i+3}u(g^{i+1}u)^{-1}$ divides both three
  and four, yielding a contradiction.
\end{proof}

The following is a generalization of Corollary~\ref{cor:C} from the
Introduction.

\begin{cor}\label{cor:C'}
  Let $G$ be a group.
\begin{enumerate}
\item Suppose that $G$ is nonabelian and that the abelianization of
  $G$ contains an element of infinite order. Then for every integer
  $M\ge 1$ there exist nontrivial elements $g_1, \dots, g_M \in G$
  such that $g_i$ is not conjugate to $g_j^{\pm 1}$ for $i \ne j$ and
  such that for every action of $G$ on an $\mathbb{R}$-tree $X$ with
  trivial arc stabilizers we have $\ell_X(g_i) = \ell_X(g_j)$.
  
\item Suppose that $G$ is nonabelian and that the abelianization of
  $G$ contains an element of order bigger than four. Then there exist
  nontrivial elements $g_1, g_2\in G$ such that $g_1$ is not conjugate
  to $g_2^{\pm 1}$ and such that for every action of $G$ on an
  $\mathbb{R}$-tree $X$ with trivial arc stabilizers we have
  $\ell_X(g_1) = \ell_X(g_2)$.
\end{enumerate}
\end{cor}
\begin{proof}
  Let $\overline{G}$ be the abelianization of $G$. For an element $x
  \in G$ we denote by $\overline{x}$ the image of $x$ in
  $\overline{G}$.

  (1) Suppose first that $G$ is nonabelian and that $\overline G$
  contains an element of infinite order. Then there exists a
  noncentral element $g$ of $G$ whose image has infinite order in
  $\overline{G}$. Indeed, suppose not. Then all noncentral elements of
  $G$ have finite order images in $\overline{G}$. Take $g_1\in G$ such
  that $\overline{g_1}$ has infinite order in $\overline G$. Then by
  assumption $g_1$ is central in $G$. Since $G$ is nonabelian, there
  exists a noncentral element $u\in G$. Again, by assumption,
  $\overline u$ has finite order. But then $g=g_1u$ is noncentral and
  has infinite order in $\overline G$, yielding a contradiction.
  
  Thus we can choose a noncentral element $g\in G$ such that
  $\overline g$ has infinite order.  Then there is some $f \in G$ such
  that $h:=f^{- 1} g^{- 1} f \ne g^{- 1}$. Choose $M \ge 1$ and put
  $g_i:= g^i h^{2M+1-i} = g^i f^{- 1} g^{i-2M-1} f$ for $i = 1, \dots,
  M$. Then $\overline{g}_i = \overline{g}^{2i - 2M-1}\ne 1$ for $1\le
  i\le M$. Since $\overline{g}$ has infinite order in $G$, for $i \ne
  j$ we have $\overline{g_i} \ne \overline{g_j}^{\pm 1}$ and therefore
  $g_i$ is not conjugate to $g_j^{\pm 1}$ in $G$. Also, since
  $\overline{g_i} \ne 1$, it follows that $g_i \ne 1$ for $1 \le i \le
  M$. Theorem~\ref{thm-pq} implies that for every action of $G$ on an
  $\mathbb{R}$-tree $X$ with trivial arc stabilizers we have
  $\ell_X(g_i) = \ell_X(g_j)$. This establishes part (1) of the
  corollary.

  (2) Suppose now that $G$ is nonabelian and that $\overline{G}$ has
  an element of order bigger than four. We may assume that
  $\overline{G}$ is torsion by part (1).
  
  We claim that in this case there exists a noncentral element $g$ of
  $G$ whose image has order bigger than four in $\overline{G}$. Let
  $g_0\in G$ be such that $\overline{g_0}$ has order bigger than four.
  If $g_0$ is noncentral, put $g=g_0$. Otherwise, choose a noncentral
  $u\in G$ and note that $\langle \overline{g_0}, \overline{u}\rangle$
  is a finite abelian group. By Lemma~\ref{ab} there is an integer $i$
  such that $\overline{g_0^iu}$ has order at least five and hence
  $g=g_0^iu$ is the desired noncentral element.

  Thus let $g\in G$ be a noncentral element such that $\overline{g}$
  has order bigger than four.  Hence there exists $f\in F$ such that
  $h:=f^{-1}g^{-1}f\ne g^{-1}$.
  
  Put $g_1=g^4h=g^4f^{- 1} g^{-1} f$, $g_2=g^3h^2=g^3f^{- 1} g^{-2}
  f$. Then $\overline{g_1}=\overline{g^3}$ and
  $\overline{g_2}=\overline{g}$. Since $\overline g$ has order bigger
  than four, the elements $\overline{g}, \overline{g}^{-1},
  \overline{g^3}$ are nontrivial and pairwise distinct in
  $\overline{G}$. Hence $g_2\ne 1$ and $g_1$ is not conjugate to
  $g_2^{\pm 1}$ in $G$.  Theorem \ref{thm-pq} again implies that for
  every action of $G$ on an $\mathbb{R}$-tree $X$ with trivial arc
  stabilizers we have $\ell_X(g_1) = \ell_X(g_2)$.
\end{proof}

\section{Trace identities}

The following statement is well-known and probably goes back to the
work of Klein in the late 19-th century (see, for
example,~\cite{Ho,Ma} for a proof).

\begin{lem}\label{klein}
  For any freely reduced word $w(a,b)\in F(a,b)$ there exists a
  polynomial $f_w\in {\mathbb Z}[x,y,z]$ such that for any field
  $\mathbb K$ and any matrices $A,B\in SL(2,\mathbb K)$

\[
tr \, w(A,B)= f_w(tr \, A, tr \, B, tr\, AB).
\]
\end{lem}

We now obtain Theorem~\ref{B} from the Introduction:

\setcounter{theore}{1}
\begin{theore}
  Let $F$ be a finitely generated free group and suppose that
  $u\equiv_c v$ in $F$. Then $u\equiv_t v$ in $F$.
\end{theore}

\begin{proof}
  Choose an arbitrary free basis $A$ of $F$. By Lemma~6.8 of
  Horowitz~\cite{Ho}, established via a careful analysis of traces
  under an explicit family of representations in $SL(2,\mathbb C)$,
  the assumption that $u\equiv_c v$ implies that $||u||_A=||v||_A$.
  Since $A$ was an arbitrary free basis of $F$, Theorem~\ref{A}
  implies that $u\equiv_{t} v$ in $F$.
\end{proof}

Lemma~6.8 of Horowitz~\cite{Ho} was strengthened by Southcott (see
Theorem~6.6 of \cite{So}) who proved that character equivalent
elements have essentially the same ``syllable structure'' with respect
to every free basis of $F$.

\begin{cor}
  Let $F$ be a finitely generated free group and let $u,v\in F$ be
  trace equivalent elements. Then for any word $w(x,y)\in F(x,y)$ we
  have $tr(w(u,v))=tr(w(v,u))$ and hence, by Theorem~\ref{B}, $w(u,v)$
  is translation equivalent to $ w(v,u)$ in $F$.
\end{cor}
\begin{proof}
  By Lemma~\ref{klein} there is a polynomial $f(r,s,t)\in \mathbb
  Z[r,s,t]$ such that for every $A,B\in SL(2,\mathbb C)$ we have $tr
  (w(A,B))=f( tr(A), tr(B), tr(AB))$.
  
  Let $\alpha: F\to SL(2,\mathbb C)$ be an arbitrary representation.
  Put $A=\alpha(u)$ and $B=\alpha(v)$. Thus $tr(A)=tr(B)$ and, of
  course, $tr(AB)=tr(BA)$. Therefore $f(tr(A), tr(B), tr(AB))=f(tr(B),
  tr(A), tr(BA)$ and hence $tr (w(A,B))= tr (w(B,A))$. Since $\alpha$
  was arbitrary, this implies that $tr(w(u,v))=tr(w(v,u))$, as
  required.

\end{proof}

\begin{prop}\label{prop:trace}
  Let $w(x,y)\in F(x,y)$ be a freely reduced word, let $\mathbb K$ be
  any field and let $A,B\in GL(2,\mathbb K)$ be arbitrary matrices.
  Then

\[
tr\, w(A,B)= tr\, w^R(A,B) \tag{$\dag$}
\]
\end{prop}

\begin{proof}
  
  We first will prove $(\dag)$ under the assumption that both $A,B\in
  SL(2,\mathbb K)$.
  
  Let $f_w(x,y,z)$ be the polynomial provided by Lemma~\ref{klein}.
  Let $A,B\in SL(2,\mathbb K)$ be arbitrary matrices. Note that
  $tr(A)=tr(A^{-1})$, $tr(B)=tr(B^{-1})$, and $tr(AB)=tr (B^{-1}
  A^{-1})=tr (A^{-1}B^{-1})$. Then by Lemma~\ref{klein}
\begin{gather*}
  tr (w(A,B))=f_w(tr(A), tr(B), tr(AB))=\\
  f_w(tr(A^{-1}), tr(B^{-1}), tr(A^{-1}B^{-1}))=tr (w(A^{-1},
  B^{-1}))= tr (w^R(A,B)),
\end{gather*}
where the last equality holds because
$w^R(a,b)=[w(a^{-1},b^{-1})]^{-1}$.

Consider now the general case of $GL(2,\mathbb K)$. Since every field
embeds in an algebraically closed field, it suffices to prove $(\dag)$
for algebraically closed field.  So let $\mathbb K$ be an
algebraically closed field, let $w\in F(x,y)$ be a freely reduced word
and let $X,Y\in GL(2,\mathbb K)$ be arbitrary. Since $\mathbb K$ is
algebraically closed, there exist nonzero $r,q\in \mathbb K$ such that
$r^2=\det(X)$ and $q^2=\det(Y)$.  Put $X_1:=X/r$ and $Y_1:=Y/q$. Then
$\det(X_1)=\det(Y_1)=1$, so that $X_1,Y_1\in SL(2,\mathbb K)$.
Therefore, by the already established fact, we have

\[
tr\, w(X_1,Y_1)=tr\, w^R(X_1,Y_1).
\]

Let $\sigma$ be the exponent sum on $x$ in $w$ (and hence in $w^R$)
and let $\tau$ be the exponent sum on $y$ in $w$ (and hence in $w^R$).
Then $w(X,Y)=w(rX_1,qY_1)=r^{\sigma}q^{\tau}w(X_1,Y_1)$ and, similarly
$w^R(X,Y)=w^R(rX_1,qY_1)=r^{\sigma}q^{\tau}w^R(X_1,Y_1)$.

Therefore

\[
tr\, w(X,Y)= r^{\sigma}q^{\tau} tr\, w(X_1,Y_1)=r^{\sigma}q^{\tau}
tr\, w^R(X_1,Y_1)=tr\, w^R(X,Y),
\]
as required.

\end{proof}

We refer the reader to~\cite{A} for a detailed discussion on the above
proposition.

\begin{rem}
  Theorem~\ref{B} and Proposition~\ref{prop:trace} immediately imply
  Theorem~\ref{C} established earlier by a direct argument.
\end{rem}

\begin{prop}\label{tr-powers}
  For any field $\mathbb K$, for any integers $p,q$ and for any
  conjugate matrices $A,B\in GL(2,\mathbb K)$ we have
\[
tr\, A^p B^q= tr \, A^qB^p.
\]
\end{prop}

\begin{proof}
  Again, we may assume that $\mathbb K$ is algebraically closed.  For
  any matrices $A,B\in SL(2,\mathbb K)$ conjugate in $GL(2,\mathbb K)$
  we have $tr\, A =tr \, B$, and hence $tr\, A^p B^q= tr \, B^pA^p= tr
  \, A^qB^p$ by Lemma~\ref{klein} applied to $w(x,y)=x^py^q$.
  
  Suppose now $A,B\in GL(2,\mathbb K)$ are conjugate. Since $\mathbb
  K$ is algebraically closed, there exists $s\in \mathbb K$ such that
  $s^2=\det(A)=\det(B)$. Then $A_1=A/s$ and $B_1=B/s$ have determinant
  $1$ and are still conjugate in $GL(2,\mathbb K)$. Therefore $tr\,
  A_1^p B_1^q= tr \, A_1^qB_1^p$. However $A^p B^q= s^{p+q} A_1^p
  B_1^q$ and $A^qB^p=s^{p+q} A_1^q B_1^p$ which implies that $tr\, A^p
  B^q= tr \, A^qB^p$, as required.
\end{proof}

\section{Tree actions in the boundary of the Teichm\"uller space}

Let $S$ be a closed surface of negative Euler characteristic, possibly
non-orientable, and let $\mathcal T(S)$ be the Teichm\"uller space of
$S$. We think of $\mathcal T(S)$ as the set of (isotopy classes of)
marked hyperbolic structures on $S$ or, equivalently, as the set of
(conjugacy classes of) free discrete and cocompact isometric actions
of $G=\pi_1(S)$ on $\mathbb H^2$. Let $\overline{\mathcal T(S)}$ be
the Thurston compactification of $\mathcal T(S)$.  The points of
$\mu\in \overline{\mathcal T(S)}-\mathcal T(S)$ are measured
laminations of $S$.  Each such lamination $\mu$ defines a small action
of $G$ on an $\mathbb R$-tree dual $X_{\mu}$ to the lift of this
lamination to $\mathbb H^2$~\cite{MS84,MS88}. We refer the reader to
{\cite{Be,Ka}} for a detailed discussion of this topic.

We next show that our results regarding two-variable palindromes also
apply to the elements of $\overline{\mathcal T(S)}$.

\begin{thm}\label{pal-t}
  Let $S$ be a closed surface (possibly non-orientable) of negative
  Euler characteristic and let $G=\pi_1(S)$. Then for any $g,h\in G$
  and for any $w(x,y)\in F(x,y)$ we have:

\begin{enumerate}
\item For each point $p\in \mathcal T(S)$, thought of as an action of
  $G$ on $\mathbb H^2$, the elements $w(g,h)$ and $w^R(g,h)$ have
  equal translation lengths.
\item For each point $\mu\in \overline{\mathcal T(S)}-\mathcal T(S)$
  we have
\[
\ell_{X_{\mu}}(w(g,h))=\ell_{X_{\mu}}(w^R(g,h)).
\]
\end{enumerate}

\end{thm}

\begin{proof}
  Let $p\in \mathcal T(S)$ and consider the corresponding action
  $\phi:G\to Isom({\mathbb H}^2)$ of $G$ by isometries on $\mathbb
  H^2$. Recall that, when $\mathbb H^2$ is considered in the
  upper-half space model, there is a canonical isometric action of
  $GL(2,\mathbb R)$ on $\mathbb H^2$ whose image is the full isometry
  group of $\mathbb H^2$.  Namely, let $A=\begin{pmatrix} a & b\\ c &
    d\end{pmatrix}\in GL(2,\mathbb R)$.  If $\det(A)>0$ then
  $Az=\frac{az+b}{cz+d}$ for any $z\in \mathbb H^2$.  If $\det(A)<0$
  then $Az=\frac{a\bar z+b}{c\bar z+d}$ for any $z\in \mathbb H^2$.
  The first case gives us an orientation-preserving isometry of
  $\mathbb H^2$ and the second case gives an orientation-reversing
  isometry.  It is well-known that for $A\in GL(2,\mathbb R)$ the
  trace $tr(A)$ and the determinant $\det(A)$ uniquely determine the
  translation length of $A$ as an isometry of the hyperbolic plane.

  Let $w(x,y)\in F(x,y)$ and let $g,h\in G$. Let $A\in GL(2,\mathbb
  R)$ represent $\phi(g)$ and let $B\in GL(2,\mathbb R)$ represent
  $\phi(h)$.  Then by Proposition~\ref{prop:trace} $tr\, w(A,B)=tr \,
  w^R(A,B)$.  Moreover, $\det w(A,B)= \det w^R(A,B)$. Therefore we
  have $\ell_{\mathbb H^2}(w(g,h))=\ell_{\mathbb H^2}(w^R(g,h))$, as
  claimed.
  
  Thus part (1) of the theorem is established.

  Now, part (1) immediately implies part (2). Indeed, recall that each
  element $p\in \mathcal T(S)$ determines a \emph{marked length
    spectrum} $\ell_p: G\to \mathbb R$ where $\ell_p(g)$ is the
  translation length of $g$ for the isometric action of $G$ on
  $\mathbb H^2$ corresponding to $p$. Then it is well-known that for
  any $\mu\in \overline{\mathcal T(S)}-\mathcal T(S)$ the
  length-function $\ell_{X_{\mu}}:G\to \mathbb R$ is projectively the
  limit of marked length spectra $\ell_{p_n}$ for some sequence of
  $p_n\in \mathcal T(S)$. That is, there exists a sequence of scalars
  $\lambda_n>0$ such that for every $f\in G$
\[
\ell_{X_{\mu}}(f)=\lim_{n\to\infty} \lambda_n \ell_{p_n}(f).
\]
Since for $f=w(g,h)$ and $f'=w^R(g,h)$ we know by (1) that
$\ell_{p_n}(f)=\ell_{p_n}(f')$ for all $n$, it follows that
$\ell_{X_{\mu}}(f)=\ell_{X_{\mu}}(f')$, as required.
\end{proof}

Together with Theorem~\ref{pal-t} the following result immediately
implies Theorem~\ref{E} from the Introduction:

\begin{thm}\label{surface}
  Let $S$ be a closed surface of negative Euler characteristic,
  possibly non-orientable, and let $G=\pi_1(S)$. Then for any
  conjugate $g,h\in G$ and for any integers $p,q$ we have:

\begin{enumerate}
\item For each point $p\in \mathcal T(S)$, thought of as an action of
  $G$ on $\mathbb H^2$, the elements $g^ph^q$ and $g^qh^p$ have equal
  translation lengths.
\item For each point $\mu\in \overline{\mathcal T(S)}-\mathcal T(S)$
  we have
\[
\ell_{X_{\mu}}(g^ph^q)=\ell_{X_{\mu}}(g^qh^p).
\]
\item If $S$ is orientable then for any conjugate elements $g,h \in
  G$, for each point $\mu\in \overline{\mathcal T(S)}-\mathcal T(S)$
  and for any positive integers $p,q,i,j$ such that $p+q=i+j$ we have
\[
\ell_{X_{\mu}}(g^ph^q)=\ell_{X_{\mu}}(g^ih^j).
\]
\end{enumerate}

\end{thm}
\begin{proof}
  Parts (1) and (2) are established exactly as Theorem~\ref{pal-t},
  but using Proposition~\ref{tr-powers} instead of
  Proposition~\ref{prop:trace}.
  
  To see that part (3) holds observe that, by well-known results, if
  $S$ is orientable then every orbit of the action of the mapping
  class group of $S$ on $\overline{\mathcal T(S)}-\mathcal T(S)$ is
  dense in $\overline{\mathcal T(S)}-\mathcal T(S)$.  In particular,
  this applies to the orbit of a point corresponding to the stable
  foliation of a pseudo-Anosov homeomorphism of $S$. Therefore the set
  of those $\mu\in \overline{\mathcal T(S)}-\mathcal T(S)$, such that
  the action of $G$ on the tree $X_{\mu}$ is free, is dense in
  $\overline{\mathcal T(S)}-\mathcal T(S)$.  Since
  Theorem~\ref{thm-pq} applies to free actions, part (3) of
  Theorem~\ref{surface} follows.
\end{proof}

\begin{rem}
  The argument used in the proof of part (3) of Theorem~\ref{surface}
  does not work in the case of non-orientable surfaces, as follows
  from the results of Danthony and Nogueira~\cite{DN}.  At the moment
  we do not know how to prove part (3) of Theorem~\ref{surface} in the
  non-orientable case.  The problem is that if $S$ is nonorientable
  and $\mu\in \overline{\mathcal T(S)}-\mathcal T(S)$ then the action
  of $G=\pi_1(S)$ on $X_{\mu}$ need not be very small. Note, however,
  that $G$ has an index two subgroup whose action on $X_{\mu}$ is very
  small.  Note also that, as Example~\ref{distinct} below shows, part
  (1) of Theorem~\ref{thm-pq} does not have a precise analogue for
  $SL_2$ trace identities. Hence it is not possible to argue as in the
  proof of Theorem~\ref{pal-t} to establish part (3) of
  Theorem~\ref{surface} for non-orientable surfaces.
\end{rem}

\section{Examples}

\begin{exmp}
  The palindrome and the $g^ph^q$ phenomena described above no longer
  hold for the class of small actions (as opposed to free or very
  small actions).

\begin{figure}[here]
  
  \psfrag{<a^2>}{$\langle a^2\rangle$} \psfrag{<a>}{$\langle
    a\rangle$} \psfrag{<1>}{$\langle 1\rangle$} \psfrag{e}{$e$}
  \psfrag{t}{$t$} \psfrag{f}{$f$} \epsfbox{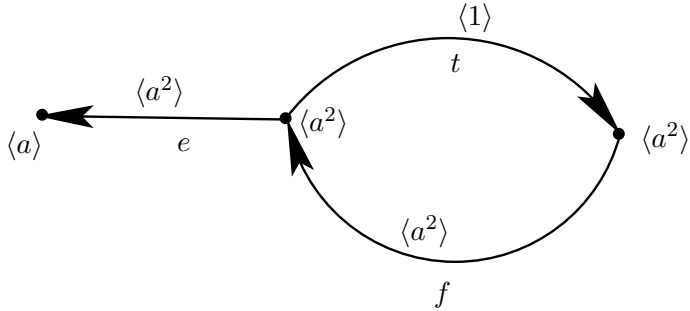}
\caption{Counterexample}\label{Fi1}

\end{figure}

Consider the graph of groups $\mathbb A$ shown in Figure~\ref{Fi1}.
Let $v$ be the vertex of the graph incident to all three marked up
edges $e,f,t$.  Put $G=\pi_1({\mathbb A},v)$ and let
$X=\widetilde{({\mathbb A},v)}$ be the Bass-Serre universal covering
tree. It is easy to see that the action of $G$ on $X$ is minimal and
that $G=F(x,y)$ is a free group of rank two with free basis $x=eae,
y=tf$.  Consider now the elements $g=xyx^2y^{-1}$ and $h=y^{-1}x^2yx$.
Thus $h$ is the palindromic reverse of $g$ in $F(x,y)$.  Moreover,
with $a=x$ and $b=yxy^{-1}$ we see that $g=ab^2$ and $h$ is conjugate
to $a^2b$.

In $G$ we have:
\begin{gather*}
  g=xyx^2y^{-1}=eae^{-1}tfea^2e^{-1}f^{-1}t^{-1}=eae^{-1}tfa^2f^{-1}t^{-1}
  =eaeta^2t^{-1}\\
  h=y^{-1}x^2yx=f^{-1}t^{-1}ea^2e^{-1}tfeae^{-1}=f^{-1}t^{-1}a^2tfeae^{-1}.
\end{gather*}
Both $eaeta^2t^{-1}$ and $f^{-1}t^{-1}a^2tfeae^{-1}$ are cyclically
reduced closed paths in $\mathbb A$.  Therefore $\ell_X(g)=4$ while
$\ell_X(h)=6$.
\end{exmp}

\begin{rem}
  We believe that the palindrome phenomenon also holds for small
  actions of a free group $F$ such that, if $g\in F$ is nontrivial and
  fixes more than one point, then $Fix(g)$ is a compact segment
  containing no branch point in its interior.
\end{rem}

\begin{exmp}
  We have remarked earlier that translation equivalence is preserved
  by injective homomorphisms, and hence by passing to larger free
  groups containing the given free group as a subgroup.  However,
  passing to subgroups, even of finite index, no longer preserves
  translation equivalence in general.
  
  Let $F=F(a,b)$ and let $H=\langle a, bab^{-1}, b^2\rangle\le F$. It
  is not hard to see that $H$ is a subgroup of index two in $F$ and
  that $x=a, y=bab^{-1}, z=b^2$ is a free basis of $H$, so that
  $H=F(x,y,z)$. Consider the elements $g=aba^2b^{-1}$ and
  $h=a^2bab^{-1}$ in $F$.  Then $g\equiv_{t} h$ in $F$ by, for
  example, Theorem~\ref{thm-pq}. We also have $g,h\in H$ and $g=xy^2,
  h=x^2y$. Clearly, $g$ and $h$ have different Whitehead graphs in
  $F(x,y,z)$ and therefore $g\not\equiv_{t} h$ in $H$.

  On the other hand, it is easy to see that if $H$ is a free factor of
  $F$ and $g,h\in H$ then $g\equiv_{t} h$ in $F$ if and only if
  $g\equiv_{t} h$ in $H$.
\end{exmp}

\begin{exmp}
  Theorem~\ref{A} states, in particular, that two elements of $F$ are
  translation equivalent in $F$ if and only if they have equal
  Whitehead graphs with respect to each free basis of $F$.  This
  statement no longer holds if we consider the obvious ``higher rank"
  analogues of Whitehead graphs where symmetrized numbers of
  occurrences of subwords of length $m > 2$ are recorded.  Again
  consider $g=aba^2b^{-1}$ and $h=a^2bab^{-1}$ in $F(a,b)$.  We
  already know that $g\equiv_{t} h$ in $F$. However, consider the
  cyclic words $w$ and $u$ defined by $g$ and $h$ accordingly. The
  number of occurrences of $b^{-1}ab$ in $w$ is equal to $1$, while
  neither $b^{-1}ab$ nor its inverse $b^{-1}a^{-1}b$ occur in $u$.
\end{exmp}

\begin{exmp}\label{distinct}
  Theorem~\ref{B} shows that in free groups character equivalence
  implies translation equivalence. However, the converse implication
  does not hold and these two phenomena are different. For example, we
  know from Theorem~\ref{D} that if $g$ and $h$ are conjugate in $F$
  and $g\ne h^{-1}$ then \newline $g^3h\equiv_{t} g^2h^2$ in $F$.  On
  the other hand, let
\[
A=\begin{bmatrix}
  2 & 1\\
  0 & 1/2
\end{bmatrix},\quad
B=\begin{bmatrix} 1 & 0\\ 2 & 1\end{bmatrix},\quad
C=BAB^{-1}=\begin{bmatrix} 0 & 1\\ -1 & 5/2\end{bmatrix}
\]
be matrices in $SL(2,\mathbb R)$.  A direct computation shows that
$tr\, A^3C=-79/16$ while $tr \, A^2C^2=-143/16$.

Thus $a^3bab^{-1}\not\equiv_c a^2 ba^2b^{-1}$ while
$a^3bab^{-1}\equiv_t a^2 ba^2b^{-1}$ in $F(a,b)$ and there is no
precise analogue of Theorem~\ref{D} for character equivalence.
\end{exmp}

\section{Open Problems}

\begin{prob}
  Is there an algorithm which, when given two elements in a finitely
  generated free group, decides whether or not they are translation
  equivalent?
\end{prob}

It can be deduced from results of Leininger~\cite{Le} that hyperbolic
equivalence in surface groups is algorithmically decidable by using
standard commutative algebra techniques applied to the representation
variety of the surface group. Similarly, one can algorithmically
decide whether or not two elements of a free group are character
equivalent.

\begin{prob}
  Find other sources of translation equivalence in free groups,
  different from those discussed in this paper.
\end{prob}

\begin{prob}
  Is it true that whenever $g\equiv_t h$ in $F$ and $w(x,y)\in F(x,y)$
  is arbitrary then $w(g,h)\equiv_t w(h,g)$ in $F$?  It easily follows
  from Lemma~\ref{klein} that $g\equiv_c h$ (e.g. $g$ is conjugate to
  $h$) implies $w(g,h)\equiv_c w(h,g)$ in $F$ and hence
  $w(g,h)\equiv_t w(h,g)$ in $F$.
\end{prob}

The notion of translation equivalence has several natural
generalizations.

\begin{prob}[Bounded translation equivalence]
  For nontrivial $g,h\in F$ we say that $g$ is \emph{boundedly
    translation equivalent} to $h$, denoted $g\equiv_b h$, if there is
  $C>0$ such that for every free and discrete action of $F$ on an
  $\mathbb R$-tree $X$ we have
\[
\frac{1}{C}\le \frac{\ell_X(g)}{\ell_X(h)}\le C.
\]

What are the sources of bounded translation equivalence in free groups
and how much more general is it compared to translation equivalence?
What are the sources of the failure of bounded translation
equivalence? Is bounded translation equivalence, in some natural
sense, generic? Is bounded translation equivalence algorithmically
decidable?
\end{prob}
A similar notion can be defined in the context of hyperbolic metrics
on closed surfaces and the above questions make sense there as well.

\begin{prob}[Volume equivalence of subgroups]
  If $G$ is a finitely generated group acting discretely isometrically
  and without a global fixed point on an $\mathbb R$-tree $X$, we
  denote by $vol_X(G)$ the sum of the lengths of the edges of the
  metric graph $X_G/G$, where $X_G$ is the unique minimal
  $G$-invariant subtree.
  
  We will say that nontrivial finitely generated subgroups $H,K$ of
  $F$ are \emph{volume equivalent} in $F$, denoted $H\equiv_v K$, if
  for every free and discrete action of $F$ on an $\mathbb R$-tree $X$
  we have $vol_X(H)=vol_X(K)$.
  
  Thus $g\equiv_t h$ in $F$ iff $\langle g\rangle \equiv_v \langle h
  \rangle$ in $F$. Moreover, if $H,K$ have the same finite index $n$
  in $F$ then it is easy to see that $H\equiv_v K$ in $F$.
  
  Are there any other sources of volume equivalence in free groups? If
  $H\equiv_v K$ in $F$, does this imply that $H$ and $K$ are free
  groups of the same rank? Is volume equivalence algorithmically
  decidable?
\end{prob}

Again, the notion of volume equivalence and the above questions also
make sense in the context of hyperbolic surfaces. In that situation
one should consider free discrete cocompact isometric actions of
$G=\pi_1(S)$ on $\mathbb H^2$. For every such action $\phi$ and for a
finitely generated non-cyclic subgroup $H$ of $G$ we define
$vol_{\phi}(H)$ as the hyperbolic volume of the compact surface
$Conv(\Lambda H)/H$ where $Conv(\Lambda H)$ is the convex hull in
$\mathbb H^2$ of the limit set $\Lambda H$ of $H$.

\end{document}